\documentclass{article}

\newtheorem{theorem}{Theorem}[section]
\newtheorem{corollary}[theorem]{Corollary}
\newtheorem{conjecture}[theorem]{Conjecture}
\newtheorem{lemma}[theorem]{Lemma}
\newtheorem{proposition}[theorem]{Proposition}

 \evensidemargin\oddsidemargin
\begin{document}

\title{Independence polynomials of well-covered graphs: generic counterexamples for
the unimodality conjecture }
\author{Vadim E. Levit and Eugen Mandrescu \\
Department of Computer Science\\
Holon Academic Institute of Technology\\
52 Golomb St., P.O. Box 305\\
Holon 58102, ISRAEL\\
\{levitv, eugen\_m\}@hait.ac.il}
\maketitle

\begin{abstract}
A graph $G$ is \textit{well-covered} if all its maximal stable sets have the
same size, denoted by $\alpha (G)$ (M. D. Plummer, 1970). If $s_{k}$ denotes
the number of stable sets of cardinality $k$ in graph $G$, and $\alpha (G)$
is the size of a maximum stable set, then $I(G;x)=\sum\limits_{k=0}^{\alpha
(G)}s_{k}x^{k}$ is the \textit{independence polynomial} of $G$ (I. Gutman
and F. Harary, 1983). J. I. Brown, K. Dilcher and R. J. Nowakowski (2000)
conjectured that $I(G;x)$ is unimodal (i.e., there is some $j\in
\{0,1,...,\alpha (G)\}$ such that $s_{0}\leq ...\leq s_{j-1}\leq s_{j}\geq
s_{j+1}\geq ...\geq s_{\alpha (G)}$) for any well-covered graph $G$. T. S.
Michael and W. N. Traves (2002) proved that this assertion is true for $%
\alpha (G)\leq 3$, while for $\alpha (G)\in \{4,5,6,7\}$ they provided
counterexamples.

In this paper we show that for any integer $\alpha \geq 8$, there
exists a connected well-covered graph $G$ with $\alpha =\alpha
(G)$, whose independence polynomial is not unimodal. In addition,
we present a number of sufficient conditions for a graph $G$ with
$\alpha (G)\leq 6$ to have the unimodal independence polynomial.

\textbf{key words:}\textit{\ stable set, independence polynomial, unimodal
sequence, well-covered graph.}
\end{abstract}

\section{Introduction}

Throughout this paper $G=(V,E)$ is a finite, undirected, loopless and
without multiple edges graph with vertex set $V=V(G)$ and edge set $E=E(G)$.
$K_{n},P_{n},K_{n_{1},n_{2},...,n_{p}}$ denote respectively, the complete
graph on $n\geq 1$ vertices, the chordless path on $n\geq 3$ vertices, and
the complete $p$-partite graph on $n_{1}+n_{2}+...+n_{p}$ vertices, where $%
n_{i}\geq 1,1\leq i\leq p$.

The \textit{disjoint union} of the graphs $G_{1},G_{2}$ is the graph $%
G=G_{1}\sqcup G_{2}$ having $V(G)=V(G_{1})\cup V(G_{2})$ and $%
E(G)=E(G_{1})\cup E(G_{2})$. In particular, $\sqcup nG$ denotes the disjoint
union of $n>1$ copies of the graph $G$. The \textit{Zykov sum} (\cite{Zykov}%
, \cite{Zykov1}) of two disjoint graphs $G_{1},G_{2}$ is the graph $%
G_{1}+G_{2}$ that has $V(G_{1})\cup V(G_{2})$ as a vertex set and $%
E(G_{1})\cup E(G_{2})\cup \{v_{1}v_{2}:v_{1}\in V(G_{1}),v_{2}\in V(G_{2})\}$
as an edge set.

A \textit{stable} set in $G$ is a set of pairwise non-adjacent vertices. The
\textit{stability number} $\alpha (G)$ of $G$ is the maximum size of a
stable set in $G$. By $\omega (G)$ we mean $\alpha (\overline{G})$, where $%
\overline{G}$ is the complement of $G$.

A graph $G$ is called \textit{well-covered} if all its maximal stable sets
are of the same cardinality, (Plummer, \cite{Plum}). If, in addition, $G$
has no isolated vertices and its order equals $2\alpha (G)$, then $G$ is
\textit{very well-covered} (Favaron, \cite{Fav1}). For instance, the graph $%
G^{*}$, obtained from $G$ by appending a single pendant edge to each vertex
of $G$ (\cite{DuttonChanBrigham}, \cite{ToppLutz}), is well-covered (see,
for example, \cite{LevMan0}), and $\alpha (G^{*})=n$. Moreover, $G^{*}$ is
very well-covered, since it is well-covered, it has no isolated vertices,
and its order equals $2\alpha (G^{*})$. The following result shows that,
under certain conditions, any well-covered graph equals $G^{*}$ for some
graph $G$.

\begin{theorem}
\cite{FinHarNow}\label{th3} Let $H$ be a connected graph of girth $\geq 6$,
which is isomorphic to neither $C_{7}$ nor $K_{1}$. Then $H$ is well-covered
if and only if its pendant edges form a perfect matching.
\end{theorem}

In other words, Theorem \ref{th3} shows that apart from $K_{1}$ and $C_{7}$,
connected well-covered graphs of girth $\geq 6$ are very well-covered. For
example, a tree $T\neq K_{1}$ could be only very well-covered, and this is
the case if and only if $T=G^{*}$ for some tree $G$ (see also Ravindra, \cite
{Ravindra}).

Let $s_{k}$ be the number of stable sets in $G$ of cardinality $k\in
\{0,1,...,\alpha (G)\}$. The polynomial $I(G;x)=\sum\limits_{k=0}^{\alpha
(G)}s_{k}x^{k}$ is called the \textit{independence polynomial} of $G$
(Gutman and Harary, \cite{GutHar}). It is easy to deduce that
\begin{eqnarray*}
I(G_{1}\sqcup G_{2};x) &=&I(G_{1};x)\cdot I(G_{2};x), \\
I(G_{1}+G_{2};x) &=&I(G_{1};x)+I(G_{2};x)-1
\end{eqnarray*}
(see also \cite{GutHar}, \cite{Arocha}, \cite{HoedeLi}).

A finite sequence of real numbers $(a_{0},a_{1},a_{2},...,a_{n})$ is said to
be:

\begin{itemize}
\item  \textit{unimodal} if there is some $k\in \{0,1,...,n\}$, called the
\textit{mode} of the sequence, such that
\[
a_{0}\leq ...\leq a_{k-1}\leq a_{k}\geq a_{k+1}\geq ...\geq a_{n},
\]

\item  \textit{log-concave} if $a_{i}^{2}\geq a_{i-1}\cdot a_{i+1}$ for $%
i\in \{1,2,...,n-1\}$.
\end{itemize}

It is known that any log-concave sequence of positive numbers is also
unimodal, while the converse is not generally true.

A polynomial $P=a_{0}+a_{1}x+a_{2}x^{2}+...+a_{n}x^{n}$ is called \textit{%
unimodal (log-concave)} if the sequence of its coefficients is unimodal
(log-concave, respectively). For instance, the independence polynomial $%
I(K_{1,3};x)=1+\mathbf{4}x+3x^{2}+x^{3}$ is log-concave, while

\[
I(K_{25}+(K_{3}\sqcup K_{4}\sqcup K_{5}\sqcup K_{5});x)=1+42x+\mathbf{107}%
x^{2}+295x^{3}+300x^{4}
\]
is unimodal, but it is not log-concave, because $107^{2}-42\cdot 295=-941$.

Hamidoune \cite{Hamidoune} proved that the independence polynomial of a
\textit{claw-free} graph (i.e., a graph having no induced subgraph
isomorphic to $K_{1,3})$ is log-concave, and hence, unimodal. However, there
are graphs whose independence polynomials are not unimodal, e.g., $%
I(K_{70}+(\sqcup 4K_{3});x)=1+\mathbf{82}x+54x^{2}+\mathbf{108}x^{3}+81x^{4}$
(for other examples, see \cite{AlMalSchErdos}). Nevertheless, in \cite
{AlMalSchErdos} it is stated the following (still open) \textit{unimodality
conjecture} for trees.

\begin{conjecture}
\label{conj2}The independence polynomial of any tree is unimodal.
\end{conjecture}

\begin{figure}[h]
\setlength{\unitlength}{1cm}%
\begin{picture}(5,2)\thicklines

  \multiput(3,0)(1,0){4}{\circle*{0.29}}
  \multiput(3,1)(1,0){4}{\circle*{0.29}}
  \multiput(3,2)(1,0){2}{\circle*{0.29}}
  \put(3,0){\line(1,0){1}}
  \put(3,1){\line(1,0){3}}
  \put(3,2){\line(1,0){1}}
  \put(4,0){\line(0,1){2}}
  \put(5,0){\line(0,1){1}}
  \put(6,0){\line(0,1){1}}

\put(2.3,1){\makebox(0,0){$T_{1}$}}

  \multiput(8,0)(1,0){3}{\circle*{0.29}}
  \multiput(8,1)(1,0){3}{\circle*{0.29}}
  \multiput(8,2)(1,0){2}{\circle*{0.29}}
  \put(8,0){\line(1,0){1}}
  \put(8,1){\line(1,0){2}}
  \put(8,2){\line(1,0){1}}
  \put(9,0){\line(0,1){2}}
  \put(10,0){\line(0,1){1}}

  \put(7.3,1){\makebox(0,0){$T_{2}$}}

 \end{picture}
\caption{Two well-covered trees.}
\label{fig5}
\end{figure}
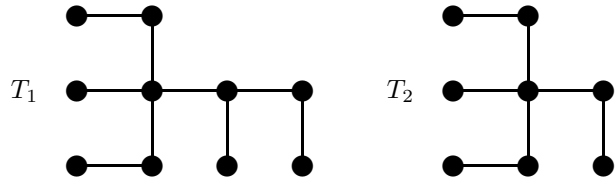

In \cite{LevitMan2} and \cite{LevitMan3}, the unimodality of independence
polynomials of a number of well-covered trees (e.g., $P_{n}^{*},K_{1,n}^{*}$%
) is validated, using the result, mentioned above, on claw-free graphs due
to Hamidoune, or directly, by identifying the location of the mode. These
findings seem promising for proving Conjecture \ref{conj2} in the case of
very well-covered trees, since a tree $T$ is well-covered if and only if
either $T$ is a \textit{well-covered spider} (i.e., $T\in
\{K_{1},K_{1}^{*},K_{1,n}^{*}:n\geq 1\}$), or $T$ is obtained from a
well-covered tree $H_{1}$ and a well-covered spider $H_{2}$, by adding an
edge joining two non-pendant vertices belonging to $H_{1},H_{2}$,\
respectively (see \cite{LevitMan1}). For instance, the trees presented in
Figure \ref{fig5} are well-covered as follows:\ $T_{2}$ is a well-covered
spider, while $T_{1}$ is an edge-join of two well-covered spiders, namely, $%
K_{1,2}^{*}$ and $K_{1,1}^{*}$.

In \cite{BrownDilNow} it was conjectured that the independence polynomial of
any well-covered graph is unimodal. Michael and Traves \cite{MichaelTraves}
proved that this assertion is true for $\alpha (G)\in \{1,2,3\}$, but it is
false for $\alpha (G)\in \{4,5,6,7\}$. Nevertheless, the conjecture of Brown
\textit{et al}. is still open for very well-covered graphs.

In \cite{LevitMan5} it was shown that for any $\alpha \geq 1$, there is a
connected very well-covered graph $G$ with $\alpha (G)=\alpha $, whose
independence polynomial is unimodal.

In this paper we prove that for any integer number $\alpha \geq 8$, there
exists a connected well-covered graph $G$ with $\alpha (G)=\alpha $, whose $%
I(G;x)$ is not unimodal. \emph{\ }We also give a simple proof for the
unimodality of the independence polynomial of a well-covered graph $G$ with $%
\alpha (G)\leq 3$, while for $\alpha (G)\in \{4,5,6\}$ a number of
sufficient conditions ensuring the unimodality of $I(G;x)$ are presented.

\section{The small stability number as a reason for well-covered graphs to
have unimodal independence polynomials}

Alavi \textit{et al.} \cite{AlMalSchErdos} showed that for any permutation $%
\sigma $ of $\{1,2,...,\alpha \}$ there is a graph $G$ with $\alpha
(G)=\alpha $ such that $s_{\sigma (1)}<s_{\sigma (2)}<...<s_{\sigma (\alpha
)}$.

\begin{lemma}
\label{lem2}If a graph $G$ satisfies $\omega (G)\leq \alpha =\alpha (G)$,
then $s_{\alpha }\leq s_{\alpha -1}$.
\end{lemma}

\setlength {\parindent}{0.0cm}\textbf{Proof. }Let $H=(\mathcal{A},\mathcal{B}%
,\mathcal{W})$ be the bipartite graph defined as follows: $X\in \mathcal{A}%
\Leftrightarrow X$ is a stable set in $G$ of size $\alpha (G)-1$, then $Y\in
\mathcal{B}\Leftrightarrow Y$ is a stable set in $G$ of size $\alpha (G)$,
and $XY\in \mathcal{W}\Leftrightarrow X\subset Y$ in $G$. Since any $Y\in
\mathcal{B}$ has exactly $\alpha (G)$ subsets of size $\alpha (G)-1$, it
follows that $\left| \mathcal{W}\right| =\alpha (G)\cdot s_{\alpha }$. On
the other hand, if $X\in \mathcal{A}$ and $X\cup \{y_{1}\},X\cup
\{y_{2}\}\in \mathcal{B}$, it implies $y_{1}y_{2}\in E(G)$, because $X$ is
stable and $\left| X\cup \{y_{1},y_{2}\}\right| >\alpha (G)$. Hence, any $%
X\in \mathcal{A}$ has at most $\omega (G)$ neighbors. Consequently, $\left|
\mathcal{W}\right| =\alpha (G)\cdot s_{\alpha }\leq \omega (G)\cdot
s_{\alpha -1}$, and this leads to $s_{\alpha }\leq s_{\alpha -1}$, since $%
\alpha (G)\geq \omega (G)$. \rule{2mm}{2mm}\setlength
{\parindent}{3.45ex}\newline

The converse of Lemma \ref{lem2} is not true, e.g., $\alpha
(K_{4}-e)=2<3=\omega (K_{4}-e)$ and $I(K_{4}-e;x)=1+4x+x^{2}$, where by $%
K_{4}-e$ we mean the graph obtained from $K_{4}$ by deleting one of its
edges.

\begin{proposition}
\cite{MichaelTraves}, \cite{LevMan4} \label{prop1}If $G$ is a well-covered
graph having $\alpha (G)=\alpha $,

then $s_{0}\leq s_{1}\leq ...\leq s_{\left\lceil \alpha /2\right\rceil }$.
\end{proposition}

\begin{corollary}
\label{cor1}If $G$ is a well-covered graph and $\omega (G)\leq \alpha (G)=3$%
, then $I(G;x)$ is log-concave.
\end{corollary}

\setlength {\parindent}{0.0cm}\textbf{Proof. }Let $%
I(G;x)=s_{0}+s_{1}x+s_{2}x^{2}+s_{3}x^{3}$. By Proposition \ref{prop1} and
Lemma \ref{lem2}, we get $s_{0}\leq s_{1}\leq s_{2}\geq s_{3}$, which
implies that $s_{2}^{2}\geq s_{1}s_{3}$. To complete the proof, let us
notice that $s_{1}^{2}=\left| V(G)\right| ^{2}\geq \left| E(\overline{G}%
)\right| =s_{2}=s_{0}s_{2}$. \rule{2mm}{2mm}\setlength
{\parindent}{3.45ex}\newline

The roots of the independence polynomials of well-covered graphs are
investigated in a number of papers, as \cite{BrownDilNow}, \cite{BrownNow},
\cite{FisherSolow}, \cite{GoldSantini}, \cite{MehHaj}, \cite{LevMan4}. Brown
\textit{et al.} showed, by a nice argument, that:

\begin{lemma}
\cite{BrownDilNow}\label{lem1} If a graph $G$ has $\alpha (G)=2$, then $%
I(G;x)$ has real roots.
\end{lemma}

The assertion fails for graphs with stability number greater than $2$, e.g.,
$I(K_{1,3};x)$. Notice that the independence polynomials of the trees from
Figure \ref{fig5}, are respectively
\begin{eqnarray*}
I(T_{1};x) &=&1+10x+36x^{2}+\mathbf{60}x^{3}+47x^{4}+14x^{5}, \\
I(T_{2};x) &=&1+8x+21x^{2}+\mathbf{23}x^{3}+9x^{4},
\end{eqnarray*}
while only for the first is true that all its roots are real. Let us observe
that $T_{1},T_{2}$ are well-covered and their polynomials are unimodal.
Hence, Newton's theorem (stating that if a polynomial with positive
coefficients has only real roots, then its coefficients form a log-concave
sequence) is not useful in solving Conjecture \ref{conj2}, even for the
particular case of very well-covered\emph{\ }trees.

Let us mention that there are connected graphs, with stability number equal
to $3$, whose independence polynomials are:

\begin{itemize}
\item  not unimodal, e.g.,
\[
I(K_{24}+(K_{4}\sqcup K_{3}\sqcup K_{3}));x)=1+\mathbf{34}x+33x^{2}+\mathbf{%
36}x^{3};
\]

\item  unimodal, but not log-concave, e.g.,
\[
I(K_{95}+(\sqcup 3K_{7}));x)=1+116x+\mathbf{147}x^{2}+343x^{3};
\]

\item  unimodal, but not log-concave, while the graphs are also
well-covered, e.g.,
\[
I((\sqcup 3K_{10})+K_{\underbrace{3,3,...,3}_{120}};x)=1+390x+%
\mathbf{660}x^{2}+1120x^{3}.
\]
\end{itemize}

There are also well-covered connected graphs with stability number equal to $%
4$, whose independence polynomials are:

\begin{itemize}
\item  not unimodal, e.g.,
\[
I((\sqcup 4K_{10})+K_{\underbrace{4,4,...,4}_{1800}};x)=1+7240x+%
\mathbf{11400}x^{2}+11200x^{3}+\mathbf{11800}x^{4};
\]

\item  unimodal, but not log-concave, e.g.,
\[
I((\sqcup 4K_{10})+K_{\underbrace{4,4,...,4}_{25}};x)=1+140x+%
\mathbf{750}x^{2}+4100x^{3}+10025x^{4};
\]

\item  log-concave, e.g.,
\[
I((\sqcup 4K_{10})+K_{\underbrace{4,4,...,4}_{10}%
};x)=1+80x+660x^{2}+4040x^{3}+10010x^{4}.
\]
\end{itemize}

Let us observe that the product of two unimodal independence polynomials is
not always unimodal, e.g., $I(K_{100}+\sqcup 3K_{7};x)=1+121x+147x^{2}+%
\mathbf{343}x^{3}$ and $I(K_{90}+\sqcup 3K_{7};x)=1+111x+147x^{2}+\mathbf{343%
}x^{3}$, while their product is not unimodal:
\[
1+232x+13725x^{2}+34790x^{3}+\mathbf{101185}x^{4}+100842x^{5}+\mathbf{117649}%
x^{6}.
\]

\begin{theorem}
\cite{KeilsonGerber}\label{th2} The product of a log-concave polynomial by a
unimodal polynomial is unimodal, while the product of two log-concave
polynomials is log-concave.
\end{theorem}

Theorem \ref{th2} is best possible for independence polynomials, since the
product of a log-concave independence polynomial and a unimodal independence
polynomial is not always log-concave. For instance, $I(K_{40}+\sqcup
3K_{7};x)=1+61x+147x^{2}+\mathbf{343}x^{3}$ is log-concave, $%
I(K_{110}+\sqcup 3K_{7};x)=1+131x+147x^{2}+\mathbf{343}x^{3}$ is unimodal,
while their product
\[
1+192x+8285x^{2}+28910x^{3}+87465x^{4}+\mathbf{100842}x^{5}+117649x^{6}
\]
is not log-concave.

Further we summarize some facts on graphs with small stability numbers.

\begin{proposition}
The following is a list of sufficient conditions ensuring that the
independence polynomial of a graph $G$ is unimodal:

\emph{(i)} any connected component $H$ of $G$ has $\alpha (H)\leq 2$;

\emph{(ii)} $\alpha (G)=3$ and $G$ is well-covered;

\emph{(iii)} $\alpha (G)=4,G$ is disconnected and well-covered;

\emph{(iv)} $\alpha (G)=5,G=H_{1}\sqcup H_{2},\alpha (H_{1})=2$ and $H_{2}$
is well-covered;

\emph{(v)} $\omega (G)\leq \alpha (G)\leq 5$ and $G$ is well-covered;

\emph{(vi) }$\alpha (G)=6,G$ is disconnected and any component $H$ of $G$
with $\alpha (H)\in \{3,4,5\}$ is well-covered and satisfies $\omega (H)\leq
\alpha (H)$.
\end{proposition}

\setlength {\parindent}{0.0cm}\textbf{Proof.} \emph{(i)} If $%
H_{1},H_{2},...,H_{k}$ are the components of $G$ and $\alpha (H_{i})\leq
2,1\leq i\leq k$, then $I(G;x)$ is unimodal, by Newton's Theorem, because $%
I(G;x)=I(H_{1};x)\cdot ...\cdot I(H_{k};x)$ and, consequently, by Lemma \ref
{lem1}, all its roots are real.\setlength
{\parindent}{3.45ex}

\emph{(ii)} If $G$ is disconnected, then $I(G;x)$ is unimodal, by part \emph{%
(i)}. Assume that $G$ is connected, and let $%
I(G;x)=1+nx+s_{2}x^{2}+s_{3}x^{3}$, where $n$ is the order of $G$. Any
vertex $v\in V(G)$ is contained in some maximum stable set of $G$, since $G$
is well-covered. Hence, $v$ has at least two neighbors in the complement $%
\overline{G}$ of $G$, which ensures that $n\leq \left| E(\overline{G}%
)\right| =s_{2}$. Consequently, $I(G;x)$ is unimodal, with the mode $2$ or $%
3 $, depending on $max\{s_{2},s_{3}\}$, respectively. Let us mention that
there are connected well-covered graphs with stability number equal to $3$,
whose independence polynomial has non-real roots, e.g., $I(K_{3,3,3};x)=%
\allowbreak 1+9x+9x^{2}+3x^{3}$ has non-real roots.

\emph{(iii)} If $G$ is disconnected and at least one of its components is a
complete graph, then $G=K_{p}\sqcup H$ and $I(G;x)=I(K_{p};x)\cdot
I(H;x)=(1+px)\cdot (1+s_{1}x+s_{2}x^{2}+s_{3}x^{3})$is unimodal, by Theorem
\ref{th2}. If none of its components is a complete graph, then $G$ has only
two components, say $H_{1}$ and $H_{2}$, and $\alpha (H_{1})=\alpha
(H_{2})=2 $. Hence, by Lemma \ref{lem1}, $I(H_{1};x),I(H_{2};x)$ have only
real roots. Therefore, $I(G;x)=I(H_{1};x)\cdot I(H_{2};x)$ is unimodal, by
Newton's Theorem.

\emph{(iv)} According to Lemma \ref{lem1} and Newton's Theorem, $I(H_{1};x)$
is log-concave. Since $G=H_{1}\sqcup H_{2}$, it follows that $%
I(G;x)=I(H_{1};x)\cdot I(H_{1};x)$. Hence, using part \emph{(ii)} and
Theorem \ref{th2}, we infer that $I(G;x)$ is unimodal.

\emph{(v) }Taking into account the parts \emph{(i),(ii)}, we may assume that
$\alpha (G)\in \{4,5\}$.

Suppose that $\alpha (G)=4$. Then, $%
I(G;x)=s_{0}+s_{1}x+s_{2}x^{2}+s_{3}x^{3}+s_{4}x^{4}$, and, according to
Proposition \ref{prop1}, we obtain that $s_{0}\leq s_{1}\leq s_{2}$, since $%
G $ is well-covered, while by Lemma \ref{lem2}, it follows that $s_{3}\geq
s_{4}$, because $\omega (G)\leq \alpha (G)$. Therefore, $I(G;x)$ is
unimodal, with the mode $2$ or $3$, depending on $\max \{s_{2},s_{3}\}$.
Now, for $\alpha (G)=5$, $%
I(G;x)=s_{0}+s_{1}x+s_{2}x^{2}+s_{3}x^{3}+s_{4}x^{4}+s_{5}x^{5}$ and
Proposition \ref{prop1} implies that $s_{0}\leq s_{1}\leq s_{2}\leq s_{3}$,
while Lemma \ref{lem2} assures that $s_{4}\geq s_{5}$, since $\alpha (G)\geq
\omega (G)$. Consequently, $I(G;x)$ is unimodal, with the mode $3$ or $4$,
depending on $\max \{s_{3},s_{4}\}$.

\emph{(vi) }If $G$ has a component $H$ with $\alpha (H)\in \{4,5\}$, this is
unique, and $\alpha (G-H)\leq 2$. Consequently, by parts \emph{(i)},\emph{(v)%
} and Theorem \ref{th2}, $I(G;x)=I(H;x)\cdot I(G-H;x)$ is unimodal. If $G$
has two components $H_{1},H_{2}$ with $\alpha (H_{1})=\alpha (H_{2})=3$,
then Corollary \ref{cor1} and Theorem \ref{th2} assure that $%
I(G;x)=I(H_{1};x)\cdot I(H_{2};x)$ is unimodal. The other cases follow
easily, by applying parts \emph{(i),(iii)} and Theorem \ref{th2}. \rule%
{2mm}{2mm}

\section{A family of well-covered graphs having non-unimodal independence
polynomials}

The independence polynomial of $H_{n}=(\sqcup
4K_{10})+K_{\underbrace{4,4,...,4}_{n}},n\geq 1$ is
\begin{eqnarray*}
I(H_{n};x) &=&n\cdot (1+x)^{4}+(1+10x)^{4}-n \\
&=&1+(40+4n)x+(600+6n)x^{2}+(4000+4n)x^{3}+(10000+n)x^{4}.
\end{eqnarray*}
Let us notice that $\alpha (H_{n})=4$ and $H_{n}$ is well-covered. Since $%
40+4n<600+6n$ is true for any $n\geq 1$, it follows that $I(H_{n};x)$ is not
unimodal whenever
\[
4000+4n<\min \{600+6n,10000+n\},
\]
which leads to $1700<n<2000$, where the case $n=1701$ is due to Michael and
Traves, \cite{MichaelTraves}. Moreover, $I(H_{n};x)$ is not log-concave only
for $23<n<2453$.

\begin{lemma}
\label{lem3}For any integer $k\geq 0$, the following polynomial is not
unimodal.
\[
\sum\limits_{i=0}^{k+4}s_{i}x^{i}=\left( 1+6844\cdot x+10806\cdot
x^{2}+10804\cdot x^{3}+11701\cdot x^{4}\right) \cdot \left( 1+1000\cdot
k\cdot x\right) ^{k}.
\]
\end{lemma}

\setlength {\parindent}{0.0cm}\textbf{Proof. }We show that $%
s_{k+2}>s_{k+3}<s_{k+4}$. Since the result is evident for $k=0$, we may
assume that $k\geq 1$.\setlength
{\parindent}{3.45ex}

Let us notice that:
\begin{eqnarray*}
s_{k+4} &=&11701\cdot 10^{3k}\cdot k^{k}, \\
s_{k+3} &=&10804\cdot 10^{3k}\cdot k^{k}+11701\cdot 10^{3(k-1)}\cdot
k^{k}=10^{3(k-1)}\cdot k^{k}\cdot 10815701, \\
s_{k+2} &=&10806\cdot 10^{3k}\cdot k^{k}+10804\cdot 10^{3(k-1)}\cdot k^{k}+
\\
&&+11701\cdot 10^{3(k-2)}\cdot k^{k-1}\cdot (k-1)\cdot 0.5 \\
&=&10^{3(k-2)}\cdot k^{k-1}\cdot (2\,16336\,19701\cdot k-11701)\cdot 0.5.
\end{eqnarray*}

Firstly, we have
\begin{eqnarray*}
s_{k+4}-s_{k+3} &=&11701\cdot 10^{3k}\cdot k^{k}-10^{3(k-1)}\cdot k^{k}\cdot
10815701 \\
&=&10^{3(k-1)}\cdot k^{k}\cdot 885299>0.
\end{eqnarray*}

Secondly, we obtain

\begin{eqnarray*}
s_{k+2}-s_{k+3} &=&10^{3(k-2)}\cdot k^{k-1}\cdot (2\,16336\,19701\cdot
k-11701)\cdot 0.5 \\
&&-10^{3(k-1)}\cdot k^{k}\cdot 10815701 \\
&=&10^{3(k-2)}\cdot k^{k-1}\cdot (2217701\cdot k-11701)\cdot 0.5>0,
\end{eqnarray*}
which completes the proof. \rule{2mm}{2mm}%
\begin{figure}[h]
\setlength{\unitlength}{1cm}%
\begin{picture}(5,8)\thicklines\put(1.5,0){\framebox(10.0,7.75){}}

  \put(6.5,6.35){\makebox(0,0){$G_{q}=\left( \sqcup qK_{1000}\right) \sqcup (\sqcup 4K_{10}+K_{\underbrace{4,4,...,4}_{1701}})$}}

   \put(6.5,5.65){\makebox(0,0){$q$ times}}
   \put(6.5,5.15) {\makebox(0,0){$\overbrace{\qquad  \qquad \qquad \qquad \qquad \qquad \qquad \qquad \qquad}$}}

   \put(3.6,4.15){\circle{3.0}}
   \put(3.6,4.15){\makebox(0,0){$K_{1000}$}}

  \multiput(5,4.15)(1,0){4}{\circle*{0.29}}
  \put(9.3,4.15){\circle{3.0}}
  \put(9.3,4.15){\makebox(0,0){$K_{1000}$}}

  \multiput(4,1.35)(5,0){2}{\circle{1.0}}
  \multiput(4,1.35)(5,0){2}{\makebox(0,0){$K_{10}$}}
  \put(4.5,1.35){\line(1,0){0.7}}
  \put(7.8,1.35){\line(1,0){0.7}}

  \put(6.5,1.35){\oval(2.6,1.6)}
  \put(6.5,1.35){\makebox(0,0){$K_{\underbrace{4,4,...,4}_{1701}}$}}

  \multiput(5.9,2.85)(1.3,0){2}{\circle{1.0}}
  \multiput(5.9,2.85)(1.3,0){2}{\makebox(0,0){$K_{10}$}}
  \put(5.9,2.13){\line(0,1){0.2}}
  \put(7.2,2.13){\line(0,1){0.2}}

 \end{picture}
\caption{Well-covered graphs with non-unimodal independence polynomials.}
\label{fig11}
\end{figure}
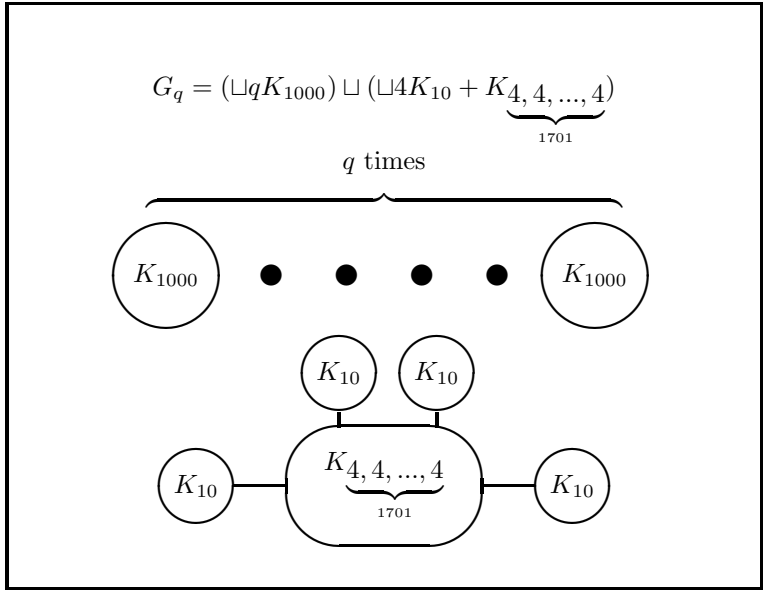

\begin{theorem}
For any integer $k\geq 4$, there is a well-covered graph $G$ with $\alpha
\left( G\right) =k$, whose independence polynomial is not unimodal.
\end{theorem}

\setlength {\parindent}{0.0cm}\textbf{Proof.} Let $q=k-4$ and $G_{q}$ be the
graph depicted in Figure \ref{fig11}, and formally defined as follows:%
\setlength
{\parindent}{3.45ex}
\[
G_{q}=\left( \sqcup qK_{1000}\right) \sqcup (\sqcup
4K_{10}+K_{\underbrace{4,4,...,4}_{1701}}).
\]

It is easy to see that $G_{q}$ is a disconnected well-covered graph, $\alpha
(G_{q})=k$, and its independence polynomial is not unimodal, because $%
I(G_{q};x)$ is identical to the non-unimodal polynomial from Lemma \ref{lem3}%
.

Moreover, the graph $G_{q}+G_{q}$ is well-covered, connected, $\alpha
(G_{q}+G_{q})=k$, and its independence polynomial is not unimodal, since $%
I(G_{q}+G_{q};x)=2\cdot I(G_{q};x)-1$. \rule{2mm}{2mm}

\section{Conclusions}

In this paper we demonstrated that for every integer $k\geq 8$ there exists
a (dis)connected well-covered graph $G$ with $\alpha (G)=k$, whose
independence polynomial is not unimodal. It is worth mentioning that all
these graphs are not very well-covered. In other words, the unimodality
conjecture remains open for the case of very well-covered graphs.

\end{document}